\documentclass[12pt]{amsart}
\usepackage{t1enc}
\usepackage[latin2]{inputenc}
\usepackage{verbatim}
\usepackage{amsmath,amsfonts,amssymb,amsthm}
\usepackage[mathcal]{eucal}
\usepackage{enumerate}
\usepackage[centertags]{amsmath}
\usepackage{graphics}

\newtheorem{theorem}{Theorem}
\newtheorem{lemma}{Lemma}

\numberwithin{equation}{subsection}

\begin{document}
\author{George Tephnadze}
\title[Fejér means ]{On the maximal operators of Walsh-Kaczmarz-Fejér means}
\address{G. Tephnadze, Department of Mathematics, Faculty of Exact and
Natural Sciences, Tbilisi State University, Chavchavadze str. 1, Tbilisi
0128, Georgia}
\email{giorgitephnadze@gmail.com}
\date{}
\maketitle

\begin{abstract}
The main aim of this paper is to prove that the maximal operator $\sigma
_{p}^{\kappa ,\ast }f:=\sup_{n\in \mathbf{P}}\left\vert \sigma _{n}^{\kappa
}f\right\vert /\left( n+1\right) ^{1/p-2}$ is bounded from the Hardy space $%
H_{p}$ to the space $L_{p},$ for $0<p<1/2$.
\end{abstract}




\section{ INTRODUCTION}

The a.e. convergence of Walsh-Fejér means $\sigma _{n}f$ was proved by Fine
\cite{F}. In 1975 Schipp \cite{Sch} showed that the maximal operator $\sigma
^{\ast }$ is of weak type $(1,1)$ and of type $(p,p)$ for $1<p\leq \infty $.
The boundedness fails to hold for $p=1$. But, Fujii \cite{Fu} proved that $%
\sigma ^{\ast }$ is bounded from the dyadic Hardy space $H_{1}$ to the space
$L_{1}$. The theorem of Fujii was extended by Weisz \cite{We4}, he showed
that the maximal operator $\sigma ^{\ast }$ is bounded from the martingale
Hardy space $H_{p}$ to the space $L_{p}$ for $p>1/2.$ Simon gave a
counterexample \cite{Si}, which shows that the boundedness does not hold for
$0<p<1/2.$ The counterexample for $p=1/2$ due to Goginava \cite{GogiEJA}
(see also \cite{tep1}). In the endpoint case $p=1/2$ two positive result was
showed. Weisz \cite{We5} proved that $\sigma ^{\ast }$ is bounded from the
Hardy space $H_{1/2}$ to the space $L_{1/2,\infty }$. In 2008 Goginava \cite%
{GoSzeged} (see also \cite{tep2}) proved that the following is true:

\textbf{TheoremA.} The maximal operator $\widetilde{\sigma }$ $^{\ast \,}$%
defined by
\begin{equation*}
\widetilde{\sigma }^{\ast }f:=\sup_{n\in \mathbf{P}}\frac{\left\vert \sigma
_{n}f\right\vert }{\log ^{2}\left( n+1\right) }
\end{equation*}%
is bounded from the Hardy space $H_{1/2}$ to the space $L_{1/2}$. He also
proved that for any nondecreasing function $\varphi :\mathbf{P}%
_{+}\rightarrow \lbrack 1,$ $\infty )$ satisfying the conditions
\begin{equation*}
\overline{\lim_{n\rightarrow \infty }}\frac{\log ^{2}\left( n+1\right) }{%
\varphi \left( n\right) }=+\infty ,
\end{equation*}%
the maximal operator
\begin{equation*}
\sup_{n\in \mathbf{P}}\frac{\left\vert \sigma _{n}f\right\vert }{\varphi
\left( n\right) }
\end{equation*}%
is not bounded from the Hardy space $H_{1/2}$ to the space $L_{1/2}.$

The author \cite{tep3} proved that when $0<p<1/2$ the maximal operator

\begin{equation}
\sigma _{p}^{\ast }f\,:=\underset{n\in \mathbf{P}}{\sup }\frac{\left\vert
\sigma _{n}f\right\vert }{\left( n+1\right) ^{1/p-2}}  \label{cond}
\end{equation}%
with respect to Walsh system is bounded from the Hardy space $H_{p}$ to the
space $L_{p}$.

We also proved that for any nondecreasing function $\varphi :\mathbf{P}%
_{+}\rightarrow \lbrack 1,$ $\infty )$ satisfying the conditions
\begin{equation}
\overline{\lim_{n\rightarrow \infty }}\frac{n^{1/p-2}}{\varphi \left(
n\right) }=+\infty ,  \label{con}
\end{equation}%
the maximal operator
\begin{equation*}
\sup_{n\in \mathbf{P}}\frac{\left\vert \sigma _{n}f\right\vert }{\varphi
\left( n\right) }
\end{equation*}%
is not bounded from the Hardy space $H_{p}$ to the space $L_{p,\infty }$
when $0<p<1/2.$ Actually, we prove a stronger result than the unboundedness
of the maximal operator $\widetilde{\sigma }_{p}^{\kappa ,\ast }f$ from the
Hardy space $H_{p}$ to the spaces $L_{p,\infty }.$ In particular, we prove
that

\begin{equation*}
\underset{n\in \mathbf{P}}{\sup }\left\Vert \frac{\sigma _{n}f}{\varphi
\left( n\right) }\right\Vert _{L_{p,\infty }}=\infty .
\end{equation*}

In 1948 $\breve{\text{S}}$neider \cite{Snei} introduced the Walsh-Kaczmarz
system and showed that the inequality
\begin{equation*}
\limsup_{n\rightarrow \infty }\frac{D_{n}^{\kappa }(x)}{\log n}\geq C>0
\end{equation*}%
holds a.e. In 1974 Schipp \cite{Sch2} and Young \cite{Y} proved that the
Walsh-Kaczmarz system is a convergence system. Skvortsov in 1981 \cite{Sk1}
showed that the Fejér means with respect to the Walsh-Kaczmarz system
converge uniformly to $f$ for any continuous functions $f$. Gát \cite{gat}
proved that, for any integrable functions, the Fejér means with respect to
the Walsh-Kaczmarz system converge almost everywhere to the function. He
showed that the maximal operator of Walsh-Kaczmarz-Fejér means $\sigma
^{\kappa ,\ast }$ is weak type $(1,1)$ and of type $(p,p)$ for all $1<p\leq
\infty $. Gát's result was generalized by Simon \cite{S2}, who showed that
the maximal operator $\sigma ^{\kappa ,\ast }$ is of type $(H_{p},L_{p})$
for $p>1/2$. In the endpoint case $p=1/2$ Goginava \cite{Gog-PM} proved that
maximal operator $\sigma ^{\kappa ,\ast }$ is not of type $(H_{1/2},L_{1/2})$
and Weisz \cite{We5} showed that the maximal operator is of weak type $%
(H_{1/2},L_{1/2})$. Analogical results of Theorem A is proved in \cite{GNCz}.

The main aim of this paper is to investigate $(H_{p},L_{p})$ and $%
(H_{p},L_{p,\infty })$-type inequalities for the maximal operators $\sigma
_{p}^{\kappa ,\ast }f:=\sup_{n\in \mathbf{P}}\left\vert \sigma _{n}^{\kappa
}f\right\vert /\left( n+1\right) ^{1/p-2}$ when $0<p<1/2.$

\section{Definitions and notations}

Now, we give a brief introduction to the theory of dyadic analysis \cite%
{SWSP}. Let $\mathbf{P}_{+}$ denote the set of positive integers, $\mathbf{%
P:=P}_{+}\cup \{0\}.$ Denote ${\mathbb{Z}}_{2}$ the discrete cyclic group of
order 2, that is ${\mathbb{Z}}_{2}=\{0,1\},$ where the group operation is
the modulo 2 addition and every subset is open. The Haar measure on ${%
\mathbb{Z}}_{2}$ is given such that the measure of a singleton is 1/2. Let $%
G $ be the complete direct product of the countable infinite copies of the
compact groups ${\mathbb{Z}}_{2}.$ The elements of $G$ are of the form $%
x=\left( x_{0},x_{1},...,x_{k},...\right) $ with $x_{k}\in \{0,1\}\left(
k\in \mathbf{P}\right) .$ The group operation on $G$ is the coordinatewise
addition, the measure (denoted by $\mu $) and the topology are the product
measure and topology. The compact Abelian group $G$ is called the Walsh
group. A base for the neighborhoods of $G$ can be given in the following
way:
\begin{equation*}
I_{0}\left( x\right) :=G,\quad I_{n}\left( x\right) :=I_{n}\left(
x_{0},...,x_{n-1}\right) :=\left\{ y\in G:\,y=\left(
x_{0},...,x_{n-1},y_{n},y_{n+1},...\right) \right\} ,
\end{equation*}
$\left( x\in G,n\in \mathbf{P}\right) .$ These sets are called dyadic
intervals. Let $0=\left( 0:i\in \mathbf{P}\right) \in G$ denote the null
element of $G,$ and $I_{n}:=I_{n}\left( 0\right) \,\left( n\in \mathbf{P}%
\right) .$ Set $e_{n}:=\left( 0,...,0,1,0,...\right) \in G,$ the $n-$th
coordinate of which is 1 and the rest are zeros $\left( n\in \mathbf{P}%
\right) .$

For $k\in \mathbf{P}$ and $x\in G$ denote
\begin{equation*}
r_{k}\left( x\right) :=\left( -1\right) ^{x_{k}}
\end{equation*}
the $k-$th Rademacher function. If $n\in \mathbf{P}$, then $%
n=\sum\limits_{i=0}^{\infty }n_{i}2^{i}$ can be written, where $n_{i}\in
\{0,1\}\quad \left( i\in \mathbf{P}\right) $, i. e. $n$ is expressed in the
number system of base 2. Denote $\left| n\right| :=\max \{j\in \mathbf{P:}%
n_{j}\neq 0\}$, that is $2^{\left| n\right| }\leq n<2^{\left| n\right| +1}.$

The Walsh-Paley system is defined as the sequence of Walsh-Paley functions:
\begin{equation*}
w_{n}\left( x\right) :=\prod\limits_{k=0}^{\infty }\left( r_{k}\left(
x\right) \right) ^{n_{k}}=r_{\left| n\right| }\left( x\right) \left(
-1\right) ^{\sum\limits_{k=0}^{\left| n\right| -1}n_{k}x_{k}},\quad \left(
x\in G,\text{ }n\in \mathbf{P}\right) .
\end{equation*}

The Walsh-Kaczmarz functions are defined by $\kappa _{0}=1$ and for $n\geq 1$

\begin{equation*}
\kappa _{n}\left( x\right) :=r_{\left| n\right| }\left( x\right)
\prod\limits_{k=0}^{\left| n\right| -1}\left( r_{\left| n\right| -1-k}\left(
x\right) \right) ^{n_{k}}=r_{\left| n\right| }\left( x\right) \left(
-1\right) ^{\sum\limits_{k=0}^{\left| n\right| -1}n_{k}x_{_{\left| n\right|
-1-k}}}.
\end{equation*}

Skvortsov (see \cite{Sk1}) gave relation between the Walsh-Kaczmarz
functions and Walsh-Paley functions by the of the transformation $\tau
_{A}:G\rightarrow G$ defined by

\begin{equation*}
\tau _{A}\left( x\right) :=\left(
x_{A-1},x_{A-2},...,x_{1},x_{0,}x_{A},x_{A+1},...\right) ,
\end{equation*}%
for $A\in \mathbf{P}.$ By the definition we have

\begin{equation*}
\kappa _{n}\left( x\right) =r_{\left| n\right| }\left( x\right)
w_{n-2^{\left| n\right| }}\left( \tau _{n}\left( x\right) \right) ,\qquad
\left( n\in \mathbf{P,}\text{ }x\in G\right) .
\end{equation*}

The Dirichlet kernels are defined

\begin{equation*}
D_{n}^{\alpha }:=\sum_{k=0}^{n-1}\alpha _{k\text{ }},
\end{equation*}
where $\alpha _{n\text{ }}=w_{n}$ or $\kappa _{n}$ $\left( n\in \mathbf{P}%
\right) ,$ $D_{0}^{\alpha }:=0.$ the $2^{n}-$th Dirichlet kernels have a
closed form (see e.g. \cite{SWSP})

\begin{equation}
D_{2^{n}}^{w}(x)=D_{2^{n}}^{\kappa }(x)=D_{2^{n}}(x)=\left\{
\begin{array}{l}
2^{n}\text{ }\,\ \ \ \,\text{if\thinspace \thinspace \thinspace }x\in I_{n}
\\
0\text{ \ \ \ \ \ \ if}\,\,x\notin I_{n}%
\end{array}%
\right.  \label{3}
\end{equation}%
The norm (or quasinorm) of the space $L_{p}(G)$ is defined by \qquad

\begin{equation*}
\left\Vert f\right\Vert _{p}:=\left( \int_{G}\left\vert f(x)\right\vert
^{p}d\mu (x)\right) ^{1/p}\qquad \left( 0<p<\infty \right) .
\end{equation*}%
The space $L_{p,\infty }(G)$ consists of all measurable functions $f$ for
which
\begin{equation*}
\left\Vert f\right\Vert _{L_{p,\infty }}:=\underset{\lambda >0}{\sup }%
\lambda \mu \left\{ \left\vert f\right\vert >\lambda \right\} ^{1/p}\leq
c<\infty .
\end{equation*}

The $\sigma $-algebra generated by the dyadic intervals of measure $2^{-k}$
will be denoted by $F_{k}$ $\left( k\in \mathbf{P}\right) .$ Denote by $%
f=\left( f^{\left( n\right) },n\in \mathbf{P}\right) $ a martingale with
respect to $\left( F_{n},n\in \mathbf{P}\right) $ (for details see, e. g.
\cite{we2}). The maximal function of a martingale $f$ is defined by
\begin{equation*}
f^{*}=\sup\limits_{n\in \mathbf{P}}\left| f^{\left( n\right) }\right| .
\end{equation*}

In case $f\in L_{1}\left( G\right) $, the maximal function can also be given
by
\begin{equation*}
f^{*}\left( x\right) =\sup\limits_{n\in \mathbf{P}}\frac{1}{\mu \left(
I_{n}(x)\right) }\left| \int\limits_{I_{n}(x)}f\left( u\right) d\mu \left(
u\right) \right| ,\ \ x\in G.
\end{equation*}

For $0<p<\infty $ the Hardy martingale space $H_{p}(G)$ consists of all
martingales for which

\begin{equation*}
\left\| f\right\| _{H_{p}}:=\left\| f^{*}\right\| _{p}<\infty .
\end{equation*}

If $f\in L_{1}\left( G\right) $, then it is easy to show that the sequence $%
\left( S_{2^{n}}f:n\in \mathbf{P}\right) $ is a martingale. If $f$ is a
martingale, that is $f=(f^{\left( 0\right) },f^{\left( 1\right) },...)$ then
the Walsh-(Kaczmarz)-Fourier coefficients must be defined in a little bit
different way:
\begin{equation*}
\widehat{f}\left( i\right) =\lim\limits_{k\rightarrow \infty
}\int\limits_{G}f^{\left( k\right) }\left( x\right) \alpha _{i}\left(
x\right) d\mu \left( x\right) ,\ \ (\alpha _{i}=w_{i}\text{ or }\kappa _{i}).
\end{equation*}

The Walsh-(Kaczmarz)-Fourier coefficients of $f\in L_{1}\left( G\right) $
are the same as the ones of the martingale $\left( S_{2^{n}}f:n\in \mathbf{P}%
\right) $ obtained from $f$.

The partial sums of the Walsh-(Kaczmarz)-Fourier series are defined as
follows:
\begin{equation*}
S_{M}^{\alpha }(f;x):=\sum\limits_{i=0}^{M-1}\widehat{f}\left( i\right)
\alpha _{i}\left( x\right) ,\quad (\alpha =w\text{ or }\kappa ).
\end{equation*}

For $n=1,2,...$ and a martingale $f$ the Fej\'er means of order $n$ of the
Walsh-(Kaczmarz)-Fourier series of the function $f$ is given by

\begin{equation*}
\sigma _{n}^{\alpha }(f;x)=\frac{1}{n}\sum\limits_{j=0}^{n-1}S_{j}^{\alpha
}(f;x),\text{ \qquad }(\alpha =w\text{ or }\kappa ).
\end{equation*}

The Fejér kernel of order $n$ of the Walsh-(Kaczmarz)-Fourier series defined
by
\begin{equation*}
K_{n}^{\alpha }\left( x\right) :=\frac{1}{n}\sum\limits_{k=0}^{n-1}D_{k}^{%
\alpha }\left( x\right) ,\text{ \qquad }(\alpha =w\text{ or }\kappa ).
\end{equation*}

For the martingale $f$ we consider maximal operators

\begin{eqnarray*}
\sigma ^{\alpha ,*}f &:&=\sup_{n\in \mathbf{P}}\left| \sigma _{n}^{\alpha
}f\right| ,\text{ \qquad }(\alpha =w\text{ or }\kappa ). \\
\overset{\sim }{\sigma }^{\alpha ,*}f &:&=\sup_{n\in \mathbf{P}}\frac{\left|
\sigma _{n}^{\alpha }f\right| }{\log ^{2}\left( n+1\right) },\text{ \qquad }%
(\alpha =w\text{ or }\kappa ). \\
\overset{\sim }{\sigma }_{p}^{\alpha ,*}f &:&=\sup_{n\in \mathbf{P}}\frac{%
\left| \sigma _{n}^{\alpha }f\right| }{\left( n+1\right) ^{1/p-2}},\text{
\qquad }(\alpha =w\text{ or }\kappa ).
\end{eqnarray*}

A bounded measurable function $a$ is p-atom, if there exists a interval I,
such that

\begin{equation*}
\left\{
\begin{array}{l}
a)\qquad \int_{I}ad\mu =0, \\
b)\ \qquad \left\| a\right\| _{\infty }\leq \mu \left( I\right) ^{\frac{-1}{p%
}}, \\
c)\qquad \text{supp}\left( a\right) \subset I.\qquad%
\end{array}
\right.
\end{equation*}

\section{FORMULATION OF MAIN RESULTS}

\begin{theorem}
Let $0<p<1/2.$ Then the \bigskip maximal operator $\overset{\sim }{\sigma }%
_{p}^{\kappa ,*}$ is bounded from the Hardy space $H_{p}\left( G\right) $ to
the space $L_{p}\left( G\right) .$
\end{theorem}

\begin{theorem}
Let $0<p<1/2$ and $\varphi :\mathbf{P}_{+}\rightarrow \lbrack 1,$ $\infty )$
be a nondecreasing function satisfying the condition
\end{theorem}

\begin{equation}
\overline{\lim_{n\rightarrow \infty }}\frac{n^{1/p-2}}{\varphi \left(
n\right) }=+\infty .  \label{6}
\end{equation}%
Then there exists a martingale $f\in H_{p}\left( G\right) ,$ such that
\begin{equation*}
\sup_{n\in \mathbf{P}}\left\Vert \frac{\sigma _{n}f}{\varphi \left( n\right)
}\right\Vert _{L_{p,\infty }\left( G\right) }=\infty .
\end{equation*}

\section{AUXILIARY PROPOSITIONS}

\begin{lemma}
(\textbf{Weisz) }\cite{we3} Suppose that an operator $T$ is sublinear and
for some $0<p\leq 1$
\end{lemma}

\begin{equation*}
\int\limits_{\overset{-}{I}}\left| Ta\right| ^{p}d\mu \leq c_{p}<\infty
\end{equation*}
for every $p$-atom $a$, where $I$ denote the support of the atom. If $T$ is
bounded from $L_{\infty \text{ }}$ to $L_{\infty \text{ }},$ then
\begin{equation*}
\left\| Tf\right\| _{L_{p}\left( G\right) }\leq c_{p}\left\| f\right\|
_{H_{p}\left( G\right) }.
\end{equation*}
\bigskip

\begin{lemma}
\textbf{(Gát) }\cite{gat} Let $A>t,$ $t,A\in N,z\in I_{t}\backslash I_{t+1}.$%
Then $\hspace*{0in}$%
\begin{equation*}
K_{2^{A}}^{w}\left( z\right) =\left\{
\begin{array}{c}
\text{ }0,\text{\qquad if }z-z_{t}e_{t}\notin I_{A}, \\
\text{ }2^{t-1},\text{\qquad if }z-z_{t}e_{t}\in I_{A}.%
\end{array}
\right.
\end{equation*}
If $z\in I_{A},$ then
\end{lemma}

\begin{equation*}
K_{2^{A}}^{w}\left( z\right) =\frac{2^{A-1}}{2}.
\end{equation*}

\begin{lemma}
\cite{GGN-SSMH} Let $n<2^{A+1},$ $A>N$ and $x\in I_{N}$ $\left(
x_{0,...,x_{m-1}},x_{m}=1,0,...,x_{l}=1,0,...,0\right) =:J_{N}^{m,l},$ $%
l=0,...,N-1,$ $m=-1,...,l.$ Then
\end{lemma}

\begin{equation*}
\int_{I_{N}}n\left\vert K_{n}^{w}\left( \tau _{A}\left( x+t\right) \right)
\right\vert d\mu \left( t\right) \leq \frac{c2^{A}}{2^{m+l}}\,,
\end{equation*}%
where $J_{N}^{-1,l}=I_{N}\left( 0,...,x_{l}=1,0,...,0\right) .$

\begin{lemma}
\textbf{(Skvortsov) }\cite{Sk1} Let $n\in \mathbf{P}_{+}\mathbf{.}$ The
following equality
\end{lemma}

\begin{eqnarray*}
nK_{n}^{\kappa }\left( x\right) &=&1+\sum\limits_{i=0}^{\left| n\right|
-1}2^{i}D_{2^{i}}\left( x\right) +\sum\limits_{i=0}^{\left| n\right|
-1}2^{i}r_{i}\left( x\right) K_{2^{i}}^{w\text{ }}\left( \tau _{i}\left(
x\right) \right) \\
&&+\left( n-2^{\left| n\right| }\right) \left( D_{2^{\left| n\right|
}}\left( x\right) +r_{\left| n\right| }\left( x\right) K_{n-2^{\left|
n\right| }}^{w\text{ }}\left( \tau _{\left| n\right| }\left( x\right)
\right) \right) .
\end{eqnarray*}
holds.

\begin{lemma}
\cite{goginava1} Let $2<A\in \mathbf{P}$ and $%
q_{A}=2^{2A}+2^{2A-2}+...+2^{2}+2^{0}.$ Then
\begin{equation*}
q_{A-1}\left\vert K_{q_{A-1}}\right\vert \geq 2^{2m+2s-3},
\end{equation*}%
for $x\in I_{2A}\left(
0,...,0,x_{2m}=1,0...,0,x_{2s}=1,x_{2s+1},...,x_{2A-1}\right) ,$ $%
m=0,1,...,A-3,$ $s=m+2,$ $m+3,...,A-1.$
\end{lemma}

\section{PROOF OF THE THEOREM}

\textbf{Proof of Theorem 1. }Lemma 4 yelds that

\begin{eqnarray*}
\overset{\sim }{\sigma }_{n}^{\kappa }f &=&\frac{f*K_{n}^{\kappa }}{\left(
n+1\right) ^{1/p-2}} \\
&\leq &\left| f*\frac{1}{\left( n+1\right) ^{1/p-1}}\left(
1+\sum\limits_{i=0}^{\left| n\right| -1}2^{i}D_{2^{i}}\right) \right| \\
&&+\left| f*\frac{1}{\left( n+1\right) ^{1/p-1}}\sum\limits_{i=0}^{\left|
n\right| -1}2^{i}r_{i}K_{2^{i}}^{w\text{ }}\circ \tau _{i}\right| \\
&&+\left| f*\frac{n-2^{\left| n\right| }}{\left( n+1\right) ^{1/p-1}}\left(
D_{2^{\left| n\right| }}+r_{\left| n\right| }K_{n-2^{\left| n\right| }}^{w%
\text{ }}\circ \tau _{\left| n\right| }\right) \right| \\
&=&f*L_{n}^{1}+f*L_{n}^{2}+f*L_{n}^{3}
\end{eqnarray*}
and

\begin{equation*}
\overset{\sim }{\sigma }_{n}^{\kappa ,*}f\leq \underset{i=1}{\overset{3}{%
\sum }}\underset{n\in P}{\sup }\left| f*L_{n}^{i}\right| =\underset{i=1}{%
\overset{3}{\sum }}R^{i}f.
\end{equation*}

The boundedness from $L_{\infty }$ to $L_{\infty }$ of the operators $R^{i}$
follows (\ref{3}) and

\begin{equation*}
\left\| K_{n}^{w}\circ \tau _{i}\right\| _{1}=\left\| K_{n}^{w}\right\|
_{1}\leq 2,
\end{equation*}
for $i<\left| n\right| ,$ $n\in \mathbf{P}$ (See \cite{Yano}). By Lemma 1
the proof will be complete, if we show that the operators $R^{i}f$ are $p$%
-quasi-local. That is, there exists a constant $c_{p}$ which is defend only $%
p$ and

\begin{equation*}
\underset{\overset{-}{I}}{\int }\left| R^{i}a\right| d\mu \leq c_{p}<\infty ,%
\text{ \qquad }i=1,2,3,
\end{equation*}
for every $p$-atom $a,$ where the dyadic interval $I$ is the support of the $%
p$-atom $a.$

Let $a$ be an arbitrary $p$-atom with support $I,$ and $\mu \left( I\right)
=2^{-N}.$ Without loss of generality, we may assume that $I:=I_{N}.$

It is evident that $\tilde{\sigma}_{n}^{\kappa }(a)=0$, if $n\leq 2^{N}.$
Therefore, we can suppose that $n>2^{N}.$

By $\Vert a\Vert _{\infty }\leq c2^{N/p}$ we have that
\begin{equation*}
|a\ast L_{n}^{i}|\leq \int_{I_{N}}|a(s)||L_{n}^{i}(x+s)|d\mu (s)\leq
c2^{N/p}\int_{I_{N}}|L_{n}^{i}(x+s)|d\mu (s)
\end{equation*}%
and
\begin{equation}
|R^{i}a|\leq c2^{N/p}\sup_{n>2^{N}}\int_{I_{N}}|L_{n}^{i}(x+s)|d\mu (s).
\label{5}
\end{equation}%
Let $x\in I_{j}\backslash I_{j+1}$ for some $j=0,...,N-1$ and $s\in I_{N},$
then $x+s\in I_{j}\backslash I_{j+1}$. Thus, we have
\begin{eqnarray*}
&&\sup_{n>2^{N}}\int_{I_{N}}|L_{n}^{1}(x+s)|d\mu (s) \\
&\leq &\sup_{n>2^{N}}\int_{I_{N}}\frac{1}{\left( n+1\right) ^{1/p-1}}\left(
1+\sum_{i=0}^{j}2^{i}D_{2^{i}}(x+s)\right) d\mu (s) \\
&\leq &\frac{c}{2^{N\left( 1/p-1\right) }}2^{2j}2^{-N}\leq \frac{c2^{2j}}{%
2^{N/p}},
\end{eqnarray*}%
Consequently,
\begin{eqnarray*}
&&\int_{\overline{I_{N}}}|R^{1}a(x)|^{p}d\mu (x) \\
&=&\sum_{j=0}^{N-1}\int_{I_{j}\backslash I_{j+1}}|R^{1}a(x)|^{p}d\mu (x)\leq
c\sum_{j=0}^{N-1}\int_{I_{j}\backslash I_{j+1}}2^{2pj}d\mu (x) \\
&\leq &c\sum_{j=0}^{N-1}2^{\left( 2p-1\right) j}\leq c<\infty ,\text{ \qquad
}0<p<1/2.
\end{eqnarray*}

Now, we discuss $\int_{\overline{I_{N}}}|R^{2}a|^{p}d\mu $. Let $\overline{%
I_{N}}=\bigcup_{t=0}^{N-1}(I_{t}\backslash I_{t+1})$ and we decompose the
sets $I_{t}\backslash I_{t+1}$ as the following disjoint union:
\begin{equation*}
I_{t}\backslash I_{t+1}=\bigcup_{l=t+1}^{N}J_{t}^{l},
\end{equation*}
where $J_{t}^{l}:=I_{N}(0,...,0,x_{t}=1,0,...,0,x_{l}=1,x_{l+1},...,x_{N-1})$
for $t<l<N$ and $J_{t}^{l}:=I_{N}(e_{t})$ for $l=N$. Then we can write
\begin{eqnarray*}
&&\int_{\overline{I_{N}}}|R^{2}a(x)|^{p}d\mu
(x)=\sum_{t=0}^{N-1}\sum_{l=t+1}^{N}\int_{J_{t}^{l}}|R^{2}a(x)|^{p}d\mu (x)
\\
&=&\sum_{t=0}^{N-1}\sum_{l=t+1}^{N-1}\int_{J_{t}^{l}}|R^{2}a(x)|^{p}d\mu
(x)+\sum_{t=0}^{N-1}\int_{J_{t}^{N}}|R^{2}a(x)|^{p}d\mu (x) \\
&=&:\sum_{1}+\sum_{2}.
\end{eqnarray*}

From Lemma 2 we have
\begin{eqnarray*}
&&\sup_{n>2^{N}}\int_{I_{N}}|L_{n}^{2}(x+s)|d\mu (s) \\
&\leq &\sup_{n>2^{N}}\frac{1}{\left( n+1\right) ^{1/p-1}}\int_{I_{N}}%
\sum_{i=0}^{l}2^{i}|K_{2^{i}}^{w}(\tau _{i}(x+s))|d\mu (s) \\
&\leq &\sup_{n>2^{N}}\frac{c}{\left( n+1\right) ^{1/p-1}}\int_{I_{N}}\left(
\sum_{i=0}^{t}2^{2i}+\sum_{i=t+1}^{l}2^{i}2^{i-t}\right) d\mu (s) \\
&\leq &\frac{c(2^{2t}+2^{2l-t})}{2^{N/p}},\text{ \qquad }x\in J_{t}^{l}.
\end{eqnarray*}

Hence,
\begin{eqnarray}
\sum_{1} &\leq &c\sum_{t=0}^{N-1}\sum_{l=t+1}^{N-1}\int_{J_{t}^{l}}\left(
2^{2t}+2^{2l-t}\right) ^{p}d\mu (x)  \label{5b} \\
&\leq &c\sum_{t=0}^{N-1}\sum_{l=t+1}^{[3t/2]}\int_{J_{t}^{l}}2^{2pt}d\mu
(x)+c\sum_{t=0}^{N-1}\sum_{l=[3t/2]+1}^{N-1}\int_{J_{t}^{l}}2^{p\left(
2l-t\right) }d\mu (x)  \notag \\
&\leq
&c\sum_{t=0}^{N-1}\sum_{l=t+1}^{[3t/2]}2^{2pt}2^{-l}+c\sum_{t=0}^{N-1}%
\sum_{l=[3t/2]+1}^{N-1}2^{p\left( 2l-t\right) }2^{-l}<c<\infty .  \notag
\end{eqnarray}

Let $x\in J_{t}^{N}$. Then Lemma 2 yields
\begin{eqnarray*}
&&\sup_{n>2^{N}}\int_{I_{N}}|L_{n}^{2}(x+s)|d\mu (s) \\
&\leq &\sup_{n>2^{N}}\frac{1}{\left( n+1\right) ^{1/p-1}}\int_{I_{N}}%
\sum_{i=0}^{|n|-1}2^{i}|K_{2^{i}}^{w}(\tau _{i}(x+s))|d\mu (s) \\
&\leq &c\sup_{n>2^{N}}\frac{1}{\left( n+1\right) ^{1/p-1}}\left(
\int_{I_{N}}\left( \sum_{i=0}^{t}2^{2i}+\sum_{i=t+1}^{N}2^{i}2^{i-t}\right)
d\mu (s)+\sum_{i=N+1}^{|n|-1}\int_{I_{i}(x_{N,i-1})}\hspace{-20pt}%
2^{i}2^{i-t}d\mu (s)\right) \\
&\leq &c\sup_{n>2^{N}}\frac{2^{2t-N}+2^{N-t}+2^{|n|-t}}{\left( n+1\right)
^{1/p-1}}\leq c\left( \frac{2^{2t}}{2^{N/p}}+\frac{2^{2N}}{2^{N/p}2^{t}}%
\right) ,
\end{eqnarray*}%
where $x_{N,i-1}:=\sum_{j=N}^{i-1}x_{j}e_{j}.$

Consequently,

\begin{eqnarray*}
\sum_{2} &\leq &c\sum_{t=0}^{N-1}\int_{J_{t}^{N}}\left( 2^{2t}+\frac{2^{2N}}{%
2^{t}}\right) ^{p}d\mu (x) \\
&\leq &c\sum_{t=0}^{[2N/3]}\int_{J_{t}^{N}}\frac{2^{2pN}}{2^{pt}}d\mu
(x)+c\sum_{t=[2N/3]+1}^{N-1}\int_{J_{t}^{N}}2^{2pt}d\mu (x) \\
&\leq &c\sum_{t=0}^{[2N/3]}\frac{1}{2^{2pt}}+c\sum_{t=[2N/3]+1}^{N-1}\frac{%
2^{2pt}}{2^{N}}\leq c.
\end{eqnarray*}

To discuss $\int_{\overline{I_{N}}}|R^{3}a|^{p}d\mu $ we use Lemma 3 and the
following disjoint decomposition of $\overline{I_{N}}:$
\begin{equation*}
\overline{I_{N}}=\bigcup_{l=0}^{N-1}\bigcup_{m=-1}^{l}J_{N}^{l,m},
\end{equation*}%
where the set $J_{N}^{l,m}$ is defined in Lemma 3.

If $x\in \overline{I_{N}}$ and $s\in I_{N}$, then $x+s\in \overline{I_{N}}$
and $D_{2^{|n|}}(x+s)=0,$ for $n>N$. Moreover, if $x\in J_{N}^{l,m}$, then $%
x+s\in J_{N}^{l,m}$ and by Lemma 3 we have
\begin{eqnarray*}
&&\sup_{n>2^{N}}\int_{I_{N}}|L_{n}^{3}(x+s)|d\mu (s) \\
&\leq &\sup_{n>2^{N}}\int_{I_{N}}\frac{n-2^{|n|}}{\left( n+1\right) ^{1/p-1}}%
|K_{n-2^{|n|}}^{w}(\tau _{|n|}(x+s))|d\mu (s) \\
&\leq &c\sup_{n>2^{N}}\frac{1}{\left( n+1\right) ^{1/p-1}}\frac{2^{|n|}}{%
2^{l+m}}\leq \frac{c2^{2N}}{2^{N/p}2^{l+m}}.
\end{eqnarray*}

Consequently,

\begin{eqnarray*}
&&\int_{\overline{I_{N}}}|R^{3}a(x)|^{p}d\mu (x) \\
&=&\sum_{l=0}^{N-1}\sum_{m=-1}^{l}\int_{J_{N}^{l,m}}|R^{3}a(x)|^{p}d\mu (x)
\\
&\leq &c\sum_{l=0}^{N-1}\sum_{m=-1}^{l}\int_{J_{N}^{l,m}}\frac{2^{2pN}}{%
2^{p(l+m)}}d\mu (x)\leq c\sum_{l=0}^{N-1}\sum_{m=-1}^{l}\frac{2^{2pN}}{%
2^{(l+m)p}}2^{-N+m} \\
&\leq &c\sum_{l=0}^{N-1}\sum_{m=-1}^{l}\frac{2^{\left( 1-p\right) m}}{%
2^{pl}2^{N\left( 1-2p\right) }}\leq c\sum_{l=0}^{N-1}\frac{2^{\left(
1-2p\right) l}}{2^{N\left( 1-2p\right) }}<c<\infty .
\end{eqnarray*}
Which complete the proof of Theorem 1.

\textbf{Proof of Theorem 2.} Let $\left\{ m_{k}:k\in \mathbf{P}\right\} $ be
an increasing sequence of positive integers such that
\begin{equation*}
\lim\limits_{k\rightarrow \infty }\frac{2^{2m_{k}\left( 1/p-2\right) }}{%
\varphi \left( q_{m_{k}}\right) }=+\infty .
\end{equation*}

Let
\begin{equation*}
f_{m_{k}}\left( x\right) :=D_{2^{2m_{k}+1}}\left( x\right)
-D_{2^{2m_{k}}}\left( x\right) .
\end{equation*}%
It is evident that
\begin{equation*}
\widehat{f}_{m_{k}}^{\kappa }\left( i\right) =\left\{
\begin{array}{l}
1,\,\,\text{if\thinspace \thinspace \thinspace }%
i=2^{2m_{k}},...,2^{2m_{k}+1}-1, \\
0,\,\,\text{otherwise.}%
\end{array}%
\right.
\end{equation*}%
Then we can write that
\begin{equation}
S_{i}^{\kappa }f_{m_{k}}\left( x\right) =\left\{
\begin{array}{ll}
D_{i}^{\kappa }\left( x\right) -D_{2^{2m_{k}}}\left( x\right) , &
i=2^{2m_{k}}+1,...,2^{2m_{k}+1}-1, \\
f_{m_{k}}\left( x\right) , & i\geq 2^{2m_{k}+1}, \\
0, & \text{otherwise.}%
\end{array}%
\right.  \label{7}
\end{equation}

Since,
\begin{equation*}
f_{m_{k}}^{\ast }\left( x\right) =\sup\limits_{n\in \mathbf{P}}\left\vert
S_{2^{n}}\left( f_{m_{k}};x\right) \right\vert =\left\vert f_{m_{k}}\left(
x\right) \right\vert ,
\end{equation*}%
from (\ref{7}) we get
\begin{equation}
\left\Vert f_{m_{k}}\right\Vert _{H_{p}}=\left\Vert f_{m_{k}}^{\ast
}\right\Vert _{p}=\left\Vert D_{2^{2m_{k}}}\right\Vert _{p}=2^{2m_{k}\left(
1-1/p\right) }.  \label{8}
\end{equation}

Since (see \cite{Sk1})

\begin{equation*}
D_{n}^{\kappa }\left( x\right) =D_{2^{\left\vert n\right\vert }}\left(
x\right) +r_{\left\vert n\right\vert }\left( x\right) D_{n-2^{\left\vert
n\right\vert }}^{w}\left( \tau _{\left\vert n\right\vert }\left( x\right)
\right)
\end{equation*}%
from (\ref{7}) we can write
\begin{eqnarray*}
&&\frac{\left\vert \sigma _{q_{m_{k}}}^{\kappa }f\left( x\right) \right\vert
}{\varphi \left( q_{m_{k}}\right) } \\
&=&\frac{1}{\varphi \left( q_{m_{k}}\right) }\frac{1}{q_{m_{k}}}\left\vert
\sum\limits_{j=0}^{q_{m_{k}}-1}S_{j}^{\kappa }f_{m_{k}}\left( x\right)
\right\vert \\
&=&\frac{1}{\varphi \left( q_{m_{k}}\right) }\frac{1}{q_{m_{k}}}\left\vert
\sum\limits_{j=2^{2m_{k}}}^{q_{m_{k}}-1}S_{j}^{\kappa }f_{m_{k}}\left(
x\right) \right\vert \\
&=&\frac{1}{\varphi \left( q_{m_{k}}\right) }\frac{1}{q_{m_{k}}}\left\vert
\sum\limits_{j=2^{2m_{k}}}^{q_{m_{k}}-1}\left( D_{j}^{\kappa }\left(
x\right) -D_{2^{2m_{k}}}\left( x\right) \right) \right\vert \\
&=&\frac{1}{\varphi \left( q_{m_{k}}\right) }\frac{1}{q_{m_{k}}}\left\vert
\sum\limits_{j=0}^{q_{m_{k}-1}-1}\left( D_{j+2^{2m_{k}}}^{\kappa }\left(
x\right) -D_{2^{2m_{k}}}\left( x\right) \right) \right\vert \\
&=&\frac{1}{\varphi \left( q_{m_{k}}\right) }\frac{1}{q_{m_{k}}}\left\vert
\sum\limits_{j=0}^{q_{m_{k}-1}-1}D_{j}^{w}\left( \tau _{2m_{k}}\left(
x\right) \right) \right\vert \\
&=&\frac{1}{\varphi \left( q_{m_{k}}\right) }\frac{q_{m_{k}-1}}{q_{m_{k}}}%
\left\vert K_{q_{m_{k}-1}}^{w}\left( \tau _{2m_{k}}\left( x\right) \right)
\right\vert
\end{eqnarray*}

Let $x\in J_{2m_{k}}^{2m_{k}-2s-1,2m_{k}-2l-1}$ for some $l<s<m_{k}.$ Then
from Lemma 5 we have

\begin{equation*}
\frac{\left\vert \sigma _{q_{m_{k}}}^{\kappa }f\left( x\right) \right\vert }{%
\varphi \left( q_{m_{k}}\right) }\geq \frac{c2^{2s+2l-2m_{k}}}{\varphi
\left( q_{m_{k}}\right) }.
\end{equation*}

Hence, we can write

\begin{eqnarray*}
&&\mu \left\{ x\in G:\frac{\left\vert \sigma _{q_{m_{k}}}^{\kappa }f\left(
x\right) \right\vert }{\varphi \left( q_{m_{k}}\right) }\geq \frac{c}{%
2^{2m_{k}}\varphi \left( q_{m_{k}}\right) }\right\} \\
&\geq &\sum\limits_{l=0}^{m_{k}-3}\sum\limits_{s=l+2}^{m_{k}-1}\mu \left\{
J_{2m_{k}}^{2m_{k}-2s-1,2m_{k}-2l-1}\right\} >c>0.
\end{eqnarray*}

Then from (\ref{8}) we obtain
\begin{eqnarray*}
&&\frac{\frac{c}{2^{2m_{k}}\varphi \left( q_{m_{k}}\right) }\left\{ \mu
\left\{ x\in G:\frac{\left\vert \sigma _{q_{m_{k}}}^{\kappa }f_{m_{k}}\left(
x\right) \right\vert }{\varphi \left( q_{m_{k}}\right) }\geq \frac{c}{%
2^{2m_{k}}\varphi \left( q_{m_{k}}\right) }\right\} \right\} ^{1/p}}{%
\left\Vert f_{m_{k}}\right\Vert _{H_{p}}} \\
&\geq &\frac{c}{2^{2m_{k}}\varphi \left( q_{m_{k}}\right) 2^{2m_{k}\left(
1-1/p\right) }} \\
&=&\frac{c2^{2m_{k}\left( 1/p-2\right) }}{\left( q_{m_{k}}\right) }%
\rightarrow \infty ,\,\,\,\,\,\text{as\thinspace \thinspace }k\rightarrow
\infty .
\end{eqnarray*}

Theorem 2 is proved.

\end{document}